\documentclass[10pt]{article}
\usepackage{latexsym}
\usepackage{graphicx}
\usepackage{color}
\usepackage{amsmath}
\usepackage{amsfonts}
\usepackage{amssymb}
\usepackage{url}
\usepackage{wrapfig}

\definecolor{darkblue}{RGB}{0, 0, 102}
\definecolor{IKB}{RGB}{0, 47, 167}

\newtheorem{theorem}{Theorem}
\newtheorem{lemma}[theorem]{Lemma}

\newtheorem{example}[theorem]{Example}

\newtheorem{remark}[theorem]{Remark}

\newtheorem{Def}[theorem]{Definition}

\newcommand{\supp}{\mbox{supp}}

\newcommand{\D}{\ensuremath{\mathcal{D}}}
\newcommand{\F}{\ensuremath{\mathcal{F}}}

\newcommand{\M}{\ensuremath{\mathcal{M}}}

\newcommand{\E}{\begin{equation}}
\newcommand{\EE}{\end{equation}}
\newcommand{\QED}{\ \rule{.1in}{.1in}}

\newcommand{\R}{\mathbb{R}}
\newcommand{\cB}{\mathbb{B}}
\newcommand{\cD}{\mathbb{D}}
\newcommand{\cF}{\mathbb{F}}
\newcommand{\cI}{\mathbb{I}}
\newcommand{\cL}{\mathbb{L}}

\newcommand{\cP}{\mathbb{P}}
\newcommand{\cR}{\mathbb{R}}
\newcommand{\cS}{\mathbb{S}}
\newcommand{\cT}{\mathbb{T}}
\newcommand{\cV}{\mathbb{V}}
\newcommand{\cX}{\mathbb{X}}
\newcommand{\cZ}{\mathbb{Z}}

\newcommand{\argmax}{\mathop{\rm argmax}}
\newcommand{\argmin}{\mathop{\rm argmin}}
\newcommand{\conv}{\mathop{\rm conv}}
\newcommand{\suppt}{\mathop{\rm suppt}}

\advance \topmargin by -\headheight
\advance \topmargin by -\headsep
\begin{document}

\Large
\begin{center}
{\bf Simpler derivation of bounded pitch inequalities for set covering, and minimum knapsack sets}\\
Daniel Bienstock and Mark Zuckerberg
\end{center}
\normalsize

\begin{abstract}
  A valid inequality $\alpha^T x \ge \alpha_0$ for a set covering problem is said to have pitch $\le \pi$ ($ \pi$ a positive integer) if the $\pi$ smallest positive $\alpha_j$ sum to at least $\alpha_0$.  This paper presents a new, simple derivation of
  a relaxation for set covering problems whose solutions satisfy all valid inequalities of pitch $\le \pi$
  and is of polynomial size, for each fixed $\pi$.  We also consider the minimum knapsack problem,
  and show that for each fixed integer $p > 0$ and $0 < \epsilon < 1$ one can separate, within additive tolerance $\epsilon$, from the relaxation defined by the
  valid inequalities with coefficients in $\{0,1,\ldots,p\}$ in time polynomial in the number of
  variables and $1/\epsilon$.
\end{abstract}

\section{Introduction}
In this paper we first consider set-covering problems.  Formally, given the $m \times n$, $0/1$ matrix $A$ we write
\begin{align}
& \Sigma(A) \ \doteq \ \{ \, x \in \{0,1\}^n \, : \, A x \ge \mbox{{\bf e}} \,\} \label{setcover}
\end{align}
\noindent where {\bf e} is the vector of $m$ 1s.  We assume, without loss of generality, that the supports
of constraints in $Ax \ge {\bf e}$  never contain one another and have cardinality greater than $1$.
Under these assumptions any undominated inequality $\alpha^T x \ge \alpha_0$ valid for $\Sigma(A)$ satisfies $\alpha_0 \ge 0$ and $\alpha_j \ge 0$ for $1 \le j \le n$.  Given such an inequality,
we write $S(\alpha) = \{ 1 \le j \le n \, : \, \alpha_j > 0\}$.

The set \eqref{setcover} can have extreme points with high fractionality that
are difficult to cut-off by standard integer programming techniques.  Balas and Ng \cite{balasng} characterize the facets for \eqref{setcover} with coefficients in $\{0,1,2\}$. The work in \cite{mark1} introduced a general class
of combinatorial inequalities generalizing those consider in \cite{balasng}.  This work relies
on the concept of \textit{pitch}:

\begin{Def}\cite{mark1}. Given an inequality $\alpha^T x \ge \alpha_0$ with $\alpha \ge 0$ and integer $0 \ge \pi \le |S(\alpha)|$, we say that $\alpha^T x \ge \alpha_0$
is of {\bf pitch} {\boldmath $\le \pi$} if the $\pi$ smallest positive
entries in $\alpha$ sum to at least $\alpha_0$. When $\pi > 0$ the inequality is said to be of pitch $= \pi$ if it is of pitch $\le \pi$ and it is not of pitch $\le \pi - 1$. 
\end{Def}
As a result, an inequality with coefficients in $\{0,1,\ldots,\pi\}$ has
pitch $\le \pi$.  Note that when $\pi = 0$ we must have $\alpha_0 \le 0$. The main result in \cite{mark1} is the following:

\begin{theorem}\label{pitched} \cite{mark1} Consider a set-covering problem given by a matrix $A$, and let $\pi \ge 2$ be a
  fixed integer.
  There is a polynomial-sized extended formulation whose projection to $x$-space  satisfies all valid inequalities for $\Sigma(A)$ of pitch $\le \pi$.\end{theorem}

Here, an \textit{extended formulation} is of the form $\{ (x,y) \in \R^{n \times N} \, : \, C x + D y \ge b \}$ (for some $N, C, D$ and $b$) whose projection to $x$-space contains $\Sigma(A)$, and the theorem states that the size of this
description (number of bits) is polynomial in $n$ and $m$. The construction in \cite{mark1} is
admittedly complex. Also see \cite{mark2}. Recent work \cite{mastro1}, \cite{fiorini1} has provided
new constructions related to bounded pitch inequalities for set covering.

A first result in the paper provides a simpler construction for Theorem \eqref{pitched}.

\begin{theorem}\label{pitchedsimple} Given a set-covering problem given by a matrix $A$, and integer $\pi \ge 2$ there is a disjunctive formulation for \eqref{setcover} with $O(m n^2 + \pi \, m^{\pi -1} n^\pi)$ variables and constraints.\end{theorem}

In the second part of this paper we consider the so-called \textit{minimum-knapsack problem},
i.e. a problem whose feasible region is of the form
\begin{align} 
  & \left\{ \, x \in \{0,1\}^n \, : \, \sum_{j = 1}^n w_j x_j \ \ge \ w_0 \, \right\} \label{minknap} 
\end{align}
with $w_j \in \cZ_+$ for $0 \le j \le n$.  Obviously by complementing the
variables we obtain a standard (``maximum'') knapsack problem; however
from a polyhedral standpoint there are significant differences.

 Without loss of generality, we assume $w_j \le w_0$ for $1 \le j \le n$.
For a vector $v \in \R^n$ and a subset of indices $S \subseteq \{1,\ldots,n\}$
we write $v(S) = \sum_{j \in S} v_j$.  As is well-known, a minimum-knapsack problem can be
equivalently restated as a set-covering problem, albeit one with an exponential number of
rows:

\begin{remark}\label{covers}Any minimum-knapsack set $\Pi \doteq \{ x \in \{0,1\}^n \, : \, \pi^T x \ge \pi_0\}$ with $\pi \ge 0$ is equivalently described by a set covering system, namely
  $$ x \in \{0,1\}^n \, : \, x(S) \ge 1 \quad \forall S \subseteq \suppt(\pi) \quad \text{with} \ \pi(S) \ge \sum_{j=1}^n \pi_j - \pi_0 + 1,$$
  i.e. the set of cover inequalities for $\Pi$.  As a corollary, an inequality
  $\alpha^T x \ge \alpha_0$ with $\alpha \ge 0$ and $\suppt(\alpha) \subseteq \suppt(\pi)$ is
  valid for $\Pi$ iff, $\forall S \subseteq \supp(\alpha)$, we have $\pi(S) \ge \sum_{j = 1}^n \pi_j - \pi_0 +1$ whenever $\alpha(S) \ge \sum_{j = 1}^n \alpha_j - \alpha_0 +1$.
\end{remark}

Here, a \textit{cover inequality} is a valid inequality of the form $x(S) \ge 1$ for some $S \subseteq \{1,\ldots,n\}$.  Given remark \eqref{covers} one wonders if a result similar to Theorem \ref{pitched} exists for the minimum-knapsack problem. This question has been taken up in recent work \cite{boundedpitchyuri}, where the following result is proved:

\begin{theorem}\label{minyuri} \cite{boundedpitchyuri} Given a minimum-knapsack problem \eqref{minknap}, $0 < \epsilon < 1$ and $p = 1$ or $p = 2$,  there is an algorithm that, with input $x^* \in [0,1]^n$  either finds a valid inequality
  for \eqref{minknap} of pitch $= p$ that is violated by $x^*$ or shows that $x^*$ satisfies
  all valid inequalities of pitch $p$ within multiplicative error $\epsilon$. The complexity of
the algorithm is polynomial in $n$ and $1/\epsilon$. \end{theorem}

This result does not guarantee strict separation; however one can show that an inequality of
pitch $= 1$ is violated iff a cover inequality is violated, and separation of such inequalities
is known to be NP-hard.  In this sense Theorem \ref{minyuri} is best possible, for pitch $\le 2$.

In this paper we prove the following result:
\begin{theorem}\label{ourminknap} Given a minimum-knapsack problem \eqref{minknap}, $0 < \epsilon < 1$, and integer $p \ge 1$,  there is an algorithm that, with input $x^* \in [0,1]^n$  either finds a valid inequality
  for \eqref{minknap} with coefficients in $\{0,1,\ldots,p\}$ that is violated by $x^*$ or shows that $x^*$ satisfies
  all such valid inequalities within additive error $\epsilon$. The complexity of
the algorithm is polynomial in $n$ and $1/\epsilon$.\end{theorem}
One can show that a valid inequality of pitch $\le 2$ is violated by a given vector $x^*$ iff
$x^*$ violates a valid inequality with coefficients in $\{0,1,2\}$; thus Theorem \ref{ourminknap}
generalizes Theorem \ref{minyuri}.

\section{Set covering}
In this section we provide a short construction leading to
Theorem \ref{pitchedsimple}.  First we will motivate our approach by introducing a branching technique for 0/1 integer programming that is of independent interest, and which we term \textit{vector branching}.

Suppose the inequality
\begin{align}
  & \sum_{j \in S} a_j x_j \ \ge \ \alpha_0 \, (> 0) \label{genM}
\end{align}
is valid for a mixed-integer set $\F$, where the $x$ variables are assumed to be binary.  Write $S \ = \ \{ j_1, j_2, \ldots, j_{|S|}\}$,
and for $1 \le h \le |S|$ define the set
$$ F_h \, \doteq \, \{ x \, \in \, \{0,1\}^{|S|} \, : \, x_{j_i} = 0 \ \forall \, i < h \mbox{ and }  x_{j_{h}} = 1\}.$$
Then the disjunction
\begin{align}
  & F_1 \, \vee \, F_2 \, \vee \ldots \vee F_{|S|} \label{vectorb}
\end{align}
is valid for $\F$ since \eqref{genM} implies $\sum_{j \in S} x_j \ge 1$.  This simple observation can be used to drive a branch-and-bound algorithm to solve optimization problems over $\F$.  When processing a node $v$ of branch-and-bound, this scheme (a) identifies an inequality \eqref{genM} that is valid for the relaxation to $\F$ used at $v$ and (b) creates $|S|$ new nodes, corresponding to each of the terms in the disjunction \eqref{vectorb}.  The node corresponding to term $h$ of \eqref{vectorb} imposes the constraints defining $F_h$ in addition to all those present at $v$. We call this procedure \textit{vector branching}\footnote{While this technique may amount to folklore, it was used in \cite{markthesis}.}.

Here we are interested in the implications of this scheme toward set covering problems \eqref{setcover} when the inequalities \eqref{genM} used to drive vector branching are rows of $Ax \ge \mbox{{\bf e}}$. Lemma \ref{vector1} given below will be used to motivate our proof of Theorem \ref{pitched}, but first,
for self-containment, we prove a basic result \cite{balasng,mark1} which will be
used in the sequel.
\begin{lemma} \label{superbasic} Suppose $\sum_{j \in T} \alpha_j x_j \ge \alpha_0$ is valid inequality for \eqref{setcover} with $\alpha \ge 0$ and $\alpha_0 > 0$. (1) There is a row  $\sum_{j \in S} \ge 1$ of $Ax \ge {\bf e}$ such that $S \subseteq T$. (2) Let $k \in T$ be such that $\alpha_k < \alpha_0$.  Then there is a row $\sum_{j \in S} \ge 1$ of $Ax \ge {\bf e}$ such that $S \subseteq T - k$.
\end{lemma}  
\noindent {\em Proof.} (1) Otherwise setting $x_j = 1$ if $j \notin T$ and $x_j = 0$ otherwise yields
a feasible binary solution to $Ax \ge {\bf e}$ which violates $\sum_{j \in T} \alpha_j x_j \ge \alpha_0$. (2) Otherwise setting $x_j = 1$ if $j \notin T$ or
$j = k$, and $x_j = 0$ otherwise, yields
a feasible binary solution to $Ax \ge {\bf e}$ which violates $\sum_{j \in T} \alpha_j x_j \ge \alpha_0$. \QED
\begin{lemma} \label{vector1} Suppose we apply vector branching to a set covering problem \eqref{setcover}. Consider a node that arises when we vector-branch on one of the rows, $\sum_{j \in S} x_j \ge 1$, of $Ax \ge {\bf e}$. Let $y \in \R^n$ be a feasible solution to the relaxation at that node. Then
  $y$ satisfies every inequality 
  \begin{align}
    & \sum_{j \in T} \alpha_j x_j \ge \alpha_0 \label{spec2}
\end{align}    
   with $\alpha_j > 0$ for all $j \in T$ which is valid for \eqref{setcover}, of pitch $\le 2$, and such that $S \subseteq T$.
\end{lemma} 
\noindent {\em Proof.} By construction $y_k = 1$ for some $k \in T$, so we may assume $\alpha_k < \alpha_0$.  By Lemma \ref{superbasic} (2) there is some row $\sum_{j \in S'} x_j \ge 1$ of $Ax \ge {\bf e}$ with $S' \subseteq T - k$.  As a result
$$ \sum_{j \in T} \alpha_j y_j \ge \sum_{j \in S'} \alpha_j y_j + \alpha_k \ge \alpha_0$$
since $\sum_{j \in S'} \alpha_j y_j \ge \min\{ \alpha_j \, : \, j \in S'\}$. \QED

This result epitomizes the techniques that we will use below to obtain a proof of Theorem \ref{pitched} -- the result relies on the explicit variable-fixing constraints defining the $F_h$ plus the structural properties of valid inequalities
for set covering.  Letchford \cite{letchford} has used a similar idea.

Note that Lemma \ref{vector1} does not yield a strategy
for developing a polynomial-size branch-and-bound tree whose leaf nodes satisfy all valid inequalities of pitch $\le 2$.  In order to prove Theorem \ref{pitched} we instead rely on the well-known equivalence between branching and disjunctive formulations.  Specifically, for each integer $\pi \ge 1$ we will
present a lifted formulation of the form
\begin{align}
  & C^\pi x \ + \ D^\pi y^\pi \ \ge \ b^\pi, \quad (x, y^\pi) \in [0,1]^n \times [0,1]^{N^\pi} \label{liftedpi}
\end{align}  
where, for some integer $N^\pi$,  $y^\pi \in $ is a vector of additional variables and $C^\pi, D^\pi, b^\pi$ are of appropriate dimension, such that
\begin{itemize}
  \item[(a)] \eqref{liftedpi} is a relaxation of \eqref{setcover}, i.e. any feasible solution to \eqref{setcover} can be lifted so as to
    satisfy \eqref{liftedpi}, 
    \item [(b)] The $x$-component of any vector $(x, y^\pi)$ feasible for \eqref{liftedpi} satisfies
      every valid inequality for \eqref{setcover} of pitch $\le \pi$, and
    \item [(c)] The size of formulation \eqref{liftedpi} is polynomial in $n$ and $m$ (= number of rows in the matrix $A$
      defining the set-covering problem).
\end{itemize}
We term
\eqref{liftedpi} the \textit{level-$\pi$ formulation}. We will
first inductively describe the level-$\pi$ formulation and then prove that it satisfies
the desired properties.  We start with the level-1 formulation, which is simply the original set-covering set of
inequalities and variable bounds: $Ax \ge {\bf e}, \ x \in [0,1]^n$ (so $N^1 = 0$).  Now assume inductively that
$\pi \ge 1$ and we have constructed the level-$\pi$ formulation. \\ 

To generate the
level-$(\pi + 1)$ formulation we proceed using steps {\bf (I)} and {\bf (II)} given below; however we first provide
an intuitive interpretation of our approach, which is based on considering the projection
to $x$-space of the level-$k$ formulation (for any $k$) which we denote by $\M^k$. 

Thus, given $\M^\pi$, for each row $i$ of $A$ we define a polytope which amounts to the disjunction
that represents vector branching on the $i^{th}$ constraint of $Ax \ge {\bf e}$.  To this end, let $1 \le i \le m$, and suppose the support of the $i^{th}$ row of $A$ is $S_i \ = \ \{j_1, \ldots, j_{|S_i|}\}$.
Then, for each $1 \le t \le |S_i|$, we define the polytope $\D^{\pi+1}(t) \subseteq [0,1]^n$ by the system
\begin{align}
  &   x_{j_t} \ =  \ 1, \quad x_{j_h} \ =  \ 0, \ \ \forall \ 1 \le h < t,  \label{projvect} \\
  &   x \ \in \ \M^{\pi}. \label{dontforget}
\end{align}
and then define
\begin{align}
  \D^{\pi+1}_i \ \doteq \ \conv \{ \D_i^{\pi+1}(t) \, : \, 1 \le t \le |S_i| \}. \label{thedisj} 
\end{align}
Finally, let
\begin{align}
  \M^{\pi+1} & \quad \doteq \quad \bigcap_i \D^{\pi+1}_i \label{key}.
\end{align}
The intuition here is that for each $1 \le i \le m$,
the collection of all systems \eqref{projvect}, \eqref{dontforget}, plus \eqref{thedisj}
indeed implements vector branching on the
$i^{th}$ constraint of $Ax \ge {\bf e}$ from a disjunctive perspective,
while \eqref{key} enforces the simultaneous application of all such disjunctions.  When $\pi = 2$, a simple
rewording of Lemma \ref{vector1} shows that every point in $\M^2$ satisfies all pitch $\le 2$ valid inequalities.\\

\noindent Now we will provide our formal description of the level-$\pi$ formulation.
\begin{itemize}
\item[{\bf (I)}] Let $1 \le i \le m$.  For $1 \le t \le |S_i|$, define the polyhedron $D^{{\pi}+1}_i(t)$ by
  the system 
\begin{subequations}  
\begin{align}
&  x_{j_t} \ =  \ 1, \quad x_{j_h} \ =  \ 0 \ \ \forall \ 1 \le h < t,  \label{zeroedYN} \\
  &  C^\pi x \, + \, D^\pi y^\pi \ \ge \ b^\pi \label{old}\\
  & (x, y^\pi) \ \in \ [0,1]^n \times [0,1]^{N^\pi}
\end{align}
\end{subequations}
Let $D_i^{\pi+1} \ \doteq \ \conv \{ D_i^{{\pi} + 1}(t) \, : \, 1 \le t \le |S_i| \}.$\\

\noindent {\bf Remark}: the projection of
$D_i^{{\pi}+1}(t)$ to $x$-space is precisely $\D^{\pi_i+1}(t)$, and likewise the projection of $D^{\pi+1}_i$ to $x$-space is
$\D^{\pi+1}_i$.\\

\noindent Formally, $D_i^{{\pi}+1}$ is described by the system
\begin{subequations}  \label{seti}
  \begin{align}
    &  \sum_{t = 1}^{|S_i|} x^{\pi, i, t}_{j_t} \ = \ 1, \quad \mbox{ and, } \label{sum}\\    
    \forall \text{ $t$ with $1 \le t \le |S_i|$:} \nonumber \\
    &  x^{\pi, i, t}_{j_h} \, - \, x^{\pi, i, t}_{j_t} \ \le \ 0 \quad \forall \, h \text{ with  } t \le h \le |S_i|, \label{x1} \\
    & x^{\pi, i, t}_{j_h} \, = \, 0 \text{ for } 1 \le h < t \label{x2}\\
&  y^{\pi, i, t}_{h} \, - \, x^{\pi, i, t}_{j_t} \ \le \ 0 \quad \forall \, h \text{ with  } 1 \le h \le N^\pi \label{y}\\
    &  C^{\pi} x^{\pi, i, t} \, + \, D^\pi y^{\pi,i,t} \ - \ x^{\pi, i, t}_{j_t} \, b^\pi \ \ge \ 0  \label{disold}\\
    & x \ = \ \sum_{t = 1}^{|S_i|} x^{\pi, i, t} \label{THEdisjunction} \\
    & x \in [0,1]^n, \ x^{\pi, i, t} \in [0,1]^n, \ y^{\pi, i, t} \in [0,1]^{N^\pi}
  \end{align}
\end{subequations}
\noindent {\bf Remark:} Constraints \eqref{x1}-\eqref{y}, and \eqref{THEdisjunction}, together with \eqref{sum} enforce the desired
vector-branching disjunction (i.e. the disjunction over all the systems \eqref{projvect}), while \eqref{disold} is the
linearization of \eqref{old}.
\item[{\bf (II)}] The level-$(\pi + 1)$ formulation is the union, over $1 \le i \le m$, of all systems \eqref{seti}.  
\end{itemize}
Next we prove the desired facts regarding this formulation.

\begin{lemma} \label{valid} Let $\hat x$ be a feasible solution to the set covering problem \eqref{setcover}, i.e. $\hat x \in \{0,1\}^n$ satisfies $A \hat x \ge {\bf e}$.  Then for every integer $\pi \ge 1$,  $\hat x$ can be lifted to a vector feasible for the level-$\pi$ formulation.
\end{lemma}
\noindent {\em Proof.} By induction on $\pi$ with the case $\pi = 1$ valid by definition. Asssume that the assertion has been proved for $\pi$; we now show it is true for $\pi + 1$.  It suffices to prove, for $1 \le i \le m$, that $\hat x$ can be lifted to a vector contained
in $D_i^{\pi + 1}$.  And in order to prove this fact we need to show that
$\hat x$ can be lifted to a vector contained in $D_i^{\pi + 1}(t)$ for some $t$ with
$1 \le t \le |S_i|$.  This fact follows by setting $t = \min \{ 1 \le h \le |S_i| \, : \, \hat x_h = 1\}$, and induction. \QED

\begin{lemma} \label{pivalid} Let $\pi \ge 1$ and suppose $(\tilde x, \tilde y)$ is a feasible solution to the level-$\pi$ formulation. Then $\hat x$ satifies   every valid inequality for \eqref{setcover} of pitch $\le \pi$.
\end{lemma}
\noindent {\em Proof.} By induction on $\pi$.  Let $\sum_{j \in T} \alpha_j x_j \ge \alpha_0$ be a valid inequality for \eqref{setcover} with $\alpha \ge 0$ and $\alpha_0 > 0$, $\alpha_j \le \alpha_0$ for each $j \in T$, of pitch $\le \pi$.  Since, for any $\pi \ge 1$, $\M^{\pi+1} \subseteq \M^{\pi}$, i.e. the set of feasible solutions for the level-$(\pi+1)$ formulation is contained in the set of feasible solutions for the level-$\pi$ formulation, we may assume that $\sum_{j \in T} \alpha_j x_j \ge \alpha_0$
has pitch = $\pi$.  Further, since $\sum_{j \in T} \alpha_j x_j \ge \alpha_0$ is valid
for \eqref{setcover}, by Lemma \ref{superbasic}(1) there exists $1 \le i \le m$ with
$S_i \subseteq T$.

Hence, if $\pi = 1$ the result clearly follows.  Suppose now that we have proved the assertion for $\pi$ and wish to prove it for $\pi+1$.  Since $(\tilde x, \tilde y)$ is contained
in $D^{\pi+1}_i$ the result will follow if we can prove that, for every $t$ with $1 \le t \le |S_i|$, any point in $D^{\pi + 1}_i(t)$ satisfies $\sum_{j \in T} \alpha_j x_j \ge \alpha_0$.

Hence, let $(\hat x, \hat y) \in D^{\pi+1}(t)$.  So $\hat x_{j_t} = 1$.  Note that $j_t \in S_i \subseteq T$. By definition, the
inequality
$$ \sum_{j \in T - j_t} \alpha_j x_j \ \ge \ \alpha_0 - \alpha_{j_t}$$
has pitch $\le \pi - 1$ and so by induction (and constraint \eqref{old}) it is satisfied
by $\hat x$.  The result now follows. \QED

\begin{lemma} \label{pisize} For each fixed $\pi$, the level-$\pi$ formulation has size
  polynomial in $m$ and $n$.
\end{lemma}
\noindent {\em Proof.} By construction of the systems $D^{\pi + 1}_i$, $N^{\pi + 1} \le m n (N^\pi + n$, which, together with $N^1 = 0$ implies $N^\pi = O(m^{\pi-1} n^{\pi-1})$.  Likewise
let $M^\pi$ denote the number of constraints in the level-$\pi$ formulation.  Then
$M^{\pi + 1} \le m ( 1 + O(n^2 + n N^\pi + n M^\pi) ) = O(m n^2  + m^{\pi -1} n^\pi + n M^\pi) = O(m n^2 + \pi \, m^{\pi -1} n^\pi)$. \QED
\section{Minimum-knapsack}
In this section we consider minimum-knapsack problems \eqref{minknap} and prove Theorem \ref{ourminknap}. 
We will use Remark \ref{covers}, together with a structural characterization of valid inequalities for \eqref{minknap}, so as to obtain a polynomial-time algorithm for near-separation from the relaxation for the knapsack \eqref{minknap} defined by the 
valid inequalities with coefficients in $\{0,1,\ldots,p\}$ for some fixed $p > 0$.
\subsubsection{Motivation: the cases $p = 2$ and $p = 3$}
We first illustrate our approach in the case $p = 2$ and outline a difficulty that arises when $p = 3$.
Let $p = 2$; here we want to check if a vector $y \in [0,1]^n$ violates any valid inequality
\begin{align}
    & x(\cS_1) + 2 x(\cS_2) \ \ge \ 2. \label{gen2}
\end{align}
(where $\cS_1 \cap \cS_2 = \emptyset$) that is not dominated by another inequality of the same kind.  Here the
situation is simple because by a result to be proven below (Theorem \ref{nodrag}) when $p =2$  we must have $w_h < w_k$ for every $h \in \cS_1$ and $k \in \cS_2$.

\begin{example}\label{2monotone}  Consider the minimum-knapsack set \eqref{minknap} given by
  inequality
\begin{align}  
 \sum_{j = 1}^n w_j x_j \ =  &\  10 x_1 + 10 x_2 + 5 x_3 + 6 x_4 + 7 x_5 \ \ge \ 10, \qquad x \in\{0,1\}^5 \nonumber
\end{align}
The inequality
\begin{align}
  &  \quad 2(x_1 + x_2 + x_3) \ + x_4 + x_5 \ \ge \ 2 \label{weak}
\end{align}
is valid.  However this inequality is not monotone in that $w_3 = 5 < 6 = w_4$ and yet the coefficients
of $x_3$ and $x_4$ are $2$ and $1$, respectively.  However, the inequality
\begin{align}  
   &  \quad 2(x_1 + x_2) \ + x_3 + x_4 + x_5 \ \ge \ 2 \nonumber 
\end{align}
is also valid, and dominates \eqref{weak}.
  \end{example}

As a result, we next argue that checking if $y$ violates any inequality \eqref{gen2} can be reduced to the solution of the following $n$ minimum-knapsack
problems, one for each index $1 \le k \le n$, where the $k^{th}$ case checks for
violations of those inequalities \eqref{gen2} where $k = \argmax\{ w_j \, : \, j \in \cS_1 \}$, and is formulated as follows:
\begin{subequations} \label{p2casek} \begin{align}
    V(k) \ &  \doteq  \ \min \quad \sum_{j \, : \, w_j \le w_k} y_j z_j \quad + \quad 2 \sum_{j \, : \, w_j > w_k} y_j z_j \label{p2obj} \\
    &  \ \text{s.t.} \qquad \qquad \sum_{j \neq k} w_j z_j \ \ge \ \sum_{j = 1}^n w_j - w_0 + 1 \label{thecover} \\
    &  \qquad \qquad z \in \{0,1\}^n,  \quad z_k = 1 \label{zconst}
\end{align} \end{subequations}
To see that this approach works, suppose $\hat z$ is feasible for \eqref{p2casek}. Then constraint \eqref{thecover} guarantees that for every index $h$ with $\hat z_h = 1$ and $w_h \le w_k$ we have
$$ \sum_{j \neq h} w_j \hat z_j \ \ge \ \sum_{j = 1}^n w_j - w_0 + 1$$
i.e. $\{j \, : \, \hat z_j = 1\}\setminus\{h\} $ forms a cover for \eqref{minknap}.  Thus, by Remark \ref{covers}, the inequality
$$ \sum_{j \, : \, w_j \le w_k, \hat z_j = 1} x_j \quad + \quad 2 \sum_{j \, : \, w_j > w_k, \hat z_j = 1} x_j  \ \ge \ 2$$
is valid for \eqref{minknap} and if it is violated by $y$ then
$$ \sum_{j \, : \, w_j \le w_k} y_j \hat z_j \quad + \quad 2 \sum_{j \, : \, w_j > w_k} y_j \hat z_j  \ < \ 2$$
and therefore $V(k) < 2$.  The reverse construction is similar. In summary, $V(k) < 2$ if and only if
$y$ violates some inequality \eqref{gen2} where $k = \argmax\{ w_j \, : \, j \in \cS_1 \}$, as desired. Note that this separation argument requires exact solution of a knapsack problem, however near-separation follows using the usual FPTAS argument.

As a further example, suppose we want to check if there exists some violated valid inequality of the form 
\begin{align}
    & x(\cT_1) + 2 x(\cT_2) \ \ge \ 3 \label{gen3}
\end{align}
with $\cT_1, \, \cT_2$ nonempty and pairwise disjoint.  [This is the $p = 3$ case but
  without a term of therm $3 x(\cT_3)$ in the left-hand side.] Again we aim to reduce this
task to a polynomially-large set of knapsack problems.
To extend the  approach used for $p = 2$ we need, to begin with, some
way to summarize the structure of covers (for \eqref{gen3}) while guaranteeing that
such covers are also covers for the original knapsack \eqref{minknap} (this will be done, in the general
case, in Lemma \ref{charqlemma} below). The salient point is that there are
two critical cases that are needed to guarantee that any cover for \eqref{gen3} is also a cover for \eqref{minknap}: first, the case where the indices of the \textit{two} largest $w_j$ with $j \in \cT_1$ are excluded from $\cT_1 \cup \cT_2$ and second, the case where
the single \textit{largest} largest $w_j$ with $j \in \cT_2$ is excluded from $\cT_1 \cup \cT_2$.  Thus, we have two cases, rather than one, but we need to be able to
represent both using a single constraint similar to \eqref{thecover}.

Moreover, the approach used for $p = 2$ relied on the ``monotonicity'' property illustrated by Example \ref{2monotone}. In the case $p = 3$ the monotonicity does not hold.
\begin{example}\label{3notmonotone}  Consider now the minimum-knapsack set \eqref{minknap} given by
\begin{align}  
 \sum_{j = 1}^n w_j x_j \ =  &\ 6 x_1 + 6 x_2 + 5 x_3 + 4 x_4 + 4 x_5 \ \ge \ 13, \qquad x \in\{0,1\}^5 \nonumber
\end{align}
The inequality
\begin{align}
  &  \quad  x_1 +  x_2 + 2 x_3 + x_4 + x_5 \ \ge \ 3 \label{undom3}
\end{align}
is valid.  Again this inequality is not monotone in that $w_3 = 5 < 6 = w_1$ and yet the coefficients
of $x_3$ and $x_1$ are $2$ and $1$, respectively.  In this case, the ``strengthening'' of \eqref{undom3} along
the lines of Example \ref{2monotone}, namely
\begin{align}  
   &  \quad x_1 + x_2  + x_3 + x_4 + x_5 \ \ge \ 2 \nonumber 
\end{align}
is \textit{not} valid.
\end{example}
As the example shows, the technique outlined for the case $p = 2$ is not easily extended due to non-monotonicities in coefficients.  However, Theorem \ref{nodrag} will show that for any fixed $p$ the total number of
such non-monotonicities is bounded (i.e. independent of $n$).

In our general approach we will handle both aspects outlined above while nevertheless relying on
enumeration of a polynomial number of cases: the multiple types of cover inequalities
that have to be considered, and the non-monotonicity illustrated by Example \ref{3notmonotone}.

\subsection{Near-monotonicity}\label{nearmonotone}
In this section we describe an important, near-monotonicity property of valid inequalities for \eqref{minknap} that will be critical in developing our polynomial-time separation algorithm.   

\begin{Def} \label{drag} Let $\sum_{j=1}^n \alpha_j x_j \ge \alpha_0$ be an inequality (valid or not) with $\alpha_j \in \cZ_+$ for $0 \le j \le n$. Let $k \in \cZ_+$ be such that $\alpha_j = k$ for some $j$. The {\bf drag} of $k$ is defined as
  $$   \ \delta(k) \ \doteq \ \left\{ h \, : \, w_h \ge \min_{j \, : \, \alpha_j = k}\{w_j\} \mbox{ and } 0 < \alpha_h < k\right\}. $$
  If there is no $j$ with $\alpha_j = k$ then we set $\delta(k) = \emptyset$. 
\end{Def}

\begin{example}\label{knapex} Consider the minimum knapsack set given by
  $$10 x_1 + 10 x_2 + 80 x_3 + 100 x_4 + 80 x_5 + 20 x_6 + 50 x_7 + 25 x_8 \ge 280, \quad x \in \{0,1\}^8.$$
  Then
  \begin{subequations}\begin{align}
      & x_1 + x_2 +  x_3 +  x_4 + 3 (x_6 +  x_7) + 4 x_8 \ge 4 \label{exq3}
  \end{align}\end{subequations}
  is valid for the knapsack.   Applying Definition \ref{drag}, some selected drag sets are as follows. $\delta(3) = \{3, 4\}$, because $w_1 = w_2 = 10 < 20 = w_6$ and both $w_3 > 20$ and $w_4 > 20$.  Similarly, $\delta(4) = \{3, 4, 7\}$.
\end{example}
Using this definition, we obtain a criterion for validity of an inequality for the knapsack set \eqref{minknap}.   
\begin{theorem} \label{nodrag}Let $P \in \cZ+$, and suppose $\sum_{j = 1}^n \alpha_j x_j \ge \alpha_0$ is a valid inequality for \eqref{minknap} with $\alpha_j \in \{0,1,\ldots,P\}$ for all $j$.  Then either $\sum_{h \in \delta(k)} \alpha_h \le P - 2$ for every $k \ge 2$ or there
  is another valid inequality with coefficients in $\{0, 1, \ldots, P\}$ that strictly dominates $\sum_{j = 1}^n \alpha_j x_j \ge \alpha_0$. \end{theorem}
\noindent {\em Proof.} Aiming for a contradiction suppose for some $k \ge 2$ we have $\sum_{h \in \delta(k)} \alpha_h  \ge P - 1$. Choose 
some index $i \in \argmin_{j \, : \, \alpha_j = k}\{w_j\}$.  For $1 \le j \le n$ write
$$ \alpha'_j  \ = \ \left\{\begin{tabular}{ll}
$\alpha_j$ & \mbox{if $j \neq i$}\\
$k - 1$ & \mbox{if $j = i$.}
\end{tabular}\right.$$
Let $\cS = \suppt(\alpha) = \suppt(\alpha')$. To complete the proof we will show that the inequality
\begin{subequations} \label{lower}\begin{align}
    & \sum_{j = 1}^n \alpha'_j x_j \ge \alpha_0,
\end{align} \end{subequations}
which dominates $\sum_{j = 1}^n \alpha_j x_j \ge \alpha_0$, is valid for \eqref{minknap}. To do so we appeal to Remark \ref{covers}.  In other words, we need to prove that
for any $C \subseteq \cS$ if $\alpha'(\cS\setminus C) < \alpha_0$ then $w(\cS\setminus C) < w_0$.

Let $C \subseteq \cS$  be given. If $i \notin \cS \setminus C$ then $\alpha'(\cS \setminus C) = \alpha(\cS \setminus C)$ and we are done. In the remainder
of the proof we assume $i \in \cS \setminus C$ and $\alpha'(\cS \setminus C) < \alpha_0$. Let $\delta = \delta(k)$ be defined as in Definition \ref{drag}.  
If $\delta \subseteq \cS \setminus C$ then 
$$ \alpha'(\cS \setminus C) \ge \alpha'(\delta) + \alpha'_i \ge P -1 + \alpha'_i \ge P $$
since $\alpha_i' = k - 1 \ge 1$.  But this is a contradiction since we assumed $\alpha'(\cS \setminus C) < \alpha_0$.  Hence $\exists h \in \delta \cap C$.
Define\footnote{For simplicity of notation in this proof we use ``+'' and ``-'' for singletons, rather than $\cap$ and $\setminus$.} $C' \doteq C + i - h$, so that $\cS \setminus C' = \cS \setminus C - i + h$. Thus $\alpha'(\cS \setminus C') = \alpha'(\cS \setminus C) - \alpha'_i + \alpha_h \le \alpha'(\cS \setminus C) < \alpha_0$ and since $\alpha(\cS \setminus C') = \alpha'(\cS \setminus C')$ we conclude that $x(C') \ge 1$ is a valid inequality for the
minimum-knapsack set defined by $\sum_{j = 1}^n \alpha_j x_j \ge \alpha_0$ and hence
it is also valid for \eqref{minknap}, i.e.
$$ w(\cS \setminus C') < w_0.$$
But $ w(\cS \setminus C) = w(\cS \setminus C') - w_h + w_i \le w(\cS \setminus C')$ since $h \in \delta$.
So $ w(\cS \setminus C) < w_0$ as desired. \QED
\subsection{Separation}  Given a
fractional vector $y \in [0,1]^n$ we consider separation of $y$ from the relaxation for the knapsack \eqref{minknap} defined by all
valid inequalities of the form
\begin{align}
    & x(\cS_1) + 2 \, x(\cS_2) + \ldots + q \, x(\cS_q) \ \ge \ q \label{genq}
\end{align}
(where the $\cS_h$ are assumed pairwise disjoint), in polynomial time. Enumeration of all $q \in \{2,\ldots, p\}$ will yield Thorem xyzp. As a first step in our procedure, we will present an efficient
procedure that succinctly enumerates all possible inequalities \eqref{genq} that are valid and
undominated.  The enumeration will be accomplished by first classifying
all inequalities \eqref{genq} using a compact scheme. Throughout, we will rely on Example \ref{knapex} given above.

\begin{example}\label{knapexb} (Example \ref{knapex}, continued.) In this example the knapsack set
  was given by the inequality $10 x_1 + 10 x_2 + 80 x_3 + 100 x_4 + 80 x_5 + 20 x_6 + 50 x_7 + 25 x_8 \ge 280$ and we consider the valid inequality $x_1 + x_2 +  x_3 +  x_4 + 3 (x_6 +  x_7) + 4 x_8 \ge 4$. Thus $q = 3$, and $\cS_1 = \{1,2,3,4\}$, $\cS_2 = \emptyset$, $\cS_3 = \{6,7\}$ and $\cS_4 = \{8\}$.
\end{example}

\noindent Consider an inequality \eqref{genq} (valid or not). In order to
classify this inequality we define
\begin{subequations}\begin{align}
    & \cI \doteq \{1 \le i \le q \ : \, S_i \neq \emptyset\}, \mbox{ and } \nonumber \\
    & \mbox{for each $i \in \cI$, } \quad \cD_i \, \doteq \, \bigcup_{q \ge k > i} \delta(k) \cap \cS_i \nonumber
\end{align}\end{subequations}    
The following observation will be useful throughout:
\begin{remark} \label{booboo2}
      Let $i \in \cI$ and suppose $j \in \cS_i \setminus \cD_i$. For any $k > i$ with $k \in \cI$ the definition of $\delta(k)$ implies $w_j < \min\{w_h \, : \, h \in \cS_k\}$. Therefore $w_j < \min\{w_h \, : \,  h \in \cD_i\}$. 
      Thus $\cD_i$ contains the indices of the $|\cD_i|$ largest $w_j$ with $j \in \cS_i$.  
\end{remark}
\noindent Now we present our classification scheme. We say that inequality \eqref{genq} (valid or not) has {\bf type} $\tau = (\cI, \cL, m)$, if $\cL = \{\cL_i \, : \, i \in \cI\}$ and  $m = \{m_i \, : \, i \in \cI \}$ satisfy:
\begin{itemize}
\item [(t.1)] For each $i \in \cI$, $\cL_i$ is a subset of $\cS_i$ satisfying
  \begin{itemize}
  \item [(a)] For any $j \in \cS_i \setminus \cL_i$ and $h \in \cL_i$ we have $w_j \le w_h$.
  \item [(b)] $| \cL_i | = \max\{ |\cD_i| \, , \, \min\{ q - 1 , |\cS_i| \} \}$.  
  \end{itemize}
      {\bf Comment}. Remark \ref{booboo2} implies  $\cD_{i} \subseteq \cL_i$. If $|\cS_i| \le q-1$ then $\cL_i = \cS_i$ and if $|\cD_i| \ge q-1$ then
      $\cL_i = \cD_i$.
 
  \noindent In Example \eqref{knapex}, $\cI = \{1, 3, 4\}$. Further $\delta(3) = \{3,4\}$, so $\cD_{1,3} = \{3, 4\}$. Similarly, $\delta(4) = \{3, 4, 7\}$, so $\cD_{1,4} = \{3, 4\}$ and $\cD_{3,4} = \{7\}$.  Hence $\cD_1 = \{3,4\}$, $\cD_3 = \{7\}$ and $\cD_4 = \emptyset$.

  Thus  $\cL_1 = \{2,3,4\}$, $\cL_3 = \{6, 7\}$ and $\cL_4 = \{8\}$ are valid choices.
\item [(t.2)] For each $k \in \cI$, $m_{k} \in \argmin\{w_j \, : \, j \in \cS_{k} \}$ with ties broken arbitrarily. 
  
  \noindent In Example \eqref{knapex},    $m_{1} = 1, m_3 = 6 \text{ and } m_4 = 8.$
\end{itemize}
\begin{remark} \label{booboo} Let $i \in \cI$.  Then 
  $\forall j \in \cS_i \setminus \cL_i, \ \ w_j \le \min\left\{ \, \min_{h \in \cL_i} w_h \ , \ \min_{i < k \in \cI} w_{m_k} - 1 \, \right\}.$
\end{remark}
  \noindent In Example \ref{knapex} we have $w_1 < w_{m_3} = w_6 = 20$, $w_1 < w_{m_4} = w_8 = 25$.\\
  \noindent Also, $\min\left\{ \, \min_{h \in \cL_1} w_h \ , \ \min_{k \in \{3,4\}} w_{m_k} - 1 \, \right\} = \min\{ 10, 20 \} = 10$.

\begin{remark} A given inequality can have more than one type. However, Lemma \ref{smallL} given next narrows the choices, without loss of generality,  and furthermore Lemma \ref{restrictsig} will show a common attribute for all types.
\end{remark}

\begin{lemma} \label{smallL} Suppose \eqref{genq} is valid for \eqref{minknap} and not dominated by
  another valid inequality with coefficients in $\{0,1,\ldots,q\}$. Then  for any type $\tau = (\cI, \cL, m)$ for \eqref{genq},  $|\cL_i| \le q^2$ for all $i \in \cI$. \end{lemma}
\noindent {\em Proof.}  Considering requirements (t.1)(a,b) for a type $\tau$, we see that for any $i \in \cI$,  $|\cL_i| \le \max \{ |\cD_i|, q - 1\}$. But $|\cD_i| \le q^2$, by Theorem \ref{nodrag}.  \QED

\noindent We will next see how the type
of an inequality encodes its validity.  This will be done in Lemma \ref{charqlemma} below after
we introduce some notation.
\begin{Def}\label{sigdef} Given an inequality \eqref{genq} of type $\tau = (\cI, \cL, m)$, its {\bf signature} is defined as 
  \begin{subequations}\label{sig}\begin{align}
      \sigma(\tau) & \doteq \ \max \ \sum_{i \in \cI} w(\cT_i) \\
      \text{s.t.} & \qquad \sum_{i \in cI} i |\cT_i| \ < \ q, \quad \text{ and } \cT_i \subseteq \cL_i, \text{ for all } i \in \cI \label{foosig}
    \end{align}
  \end{subequations}
\end{Def}
\begin{example}\label{knapex3} Example \ref{knapex}, continued. 
  In constraint \eqref{foosig} we must have $T_4 = \emptyset$, i.e. the constraint
  reads $|\cT_1| + 3|\cT_3| \le 3$ and so either $|\cT_1| = 0$ and $|\cT_3| = 1$ or
  $|\cT_1| \le 3$ and $|\cT_3| = 0$.  Clearly we obtain $\sigma = 190$.
\end{example}

\begin{lemma} \label{restrictsig} Suppose an inequality \eqref{genq} is of type $\tau$,  and that for each $i \in \cI$ we have a subset $X_i$ such that $\cL_i \subseteq \cX_i \subseteq \cS_i$.  Then we can rewrite
  \begin{subequations}\label{sigsig}\begin{align}
      \sigma(\tau) & = \ \max \ \sum_{i \in \cI} w(\cT_i) \\
      \text{s.t.} & \qquad \sum_{i \in \cI} i |\cT_i| \ < \ q, \quad \text{ and } \cT_i \subseteq \cX_i, \text{ for all } i \in \cI 
  \end{align}\end{subequations}
As a
  corollary, all types for a given inequality \eqref{genq} have the same signature.  
\end{lemma}
  \noindent {\em Proof.} Identity \eqref{sigsig} follows because  as noted in Remark \ref{booboo}
  for all $ i \in \cI$, $\cL_i$ contains the indices of the
  $| \cL_i|$-largest $w_h$ with $h \in \cS_{i}$, and by (t.1)(b) $|\cL_i| \ge \min\{q-1, |\cS_{i}|\} \ge \min\{q-1, |\cX_{i}|\}$. The corollary follows (for example) by setting $\cX_i = \cS_i$ for all $i \in \cI$.  \QED\\

\noindent We now present the characterization of validity that we will useful below.

\begin{lemma} \label{charqlemma} An inequality \eqref{genq} of type $\tau$ is valid for \eqref{minknap} iff 
\begin{align}
    \sum_{i = 1}^q w(\cS_i)  \quad & \ge \quad \sigma(\tau) \ + \ \sum_{j = 1}^n w_j \ - \ w_0 \quad + 1 \label{charq}
\end{align}
\end{lemma}
\noindent \textit{Proof.} Remark \eqref{covers} implies that inequality \eqref{genq} is valid for \eqref{minknap} iff for each family of subsets
$\cT_i \subseteq \cS_i$ ($1 \le i \le q$) such
that $\sum_{i = 1}^q i |\cT_i| < q$ we have
$$ \sum_{i = 1}^q w(\cS_i \setminus \cT_i) \ \ge \ \sum_{j = 1}^n w_j \ - \ w_0 \quad + 1 $$
from which the result follows. \QED
\begin{example}\label{knapex2} Consider Example \ref{knapex}. Inequality \eqref{exq3} is
  valid because $w(\cS_1 \cup \cS_2 \cup \cS_3 \cup \cS_4) = 295$ while
  $\sigma(\tau) + \sum_{j = 1}^8 w_j - w_0 + 1 \, = \, 190 + 375 - 280 + 1 \, = \, 286$, i.e. condition \eqref{charq} is verified.
\end{example}

\subsubsection{Separation through type enumeration}

Our separation procedure will enumerate a set of candidate triple $(\cI, \cL, m)$ that includes all possible types $\tau$ arising from valid inequalities, and for each enumerated candidate perform a
polynomial-time test, given in Section \ref{sepgiventype}.   In this section we describe the enumeration.  Let us consider an arbitrary triple $(\cI, \cL, m)$ where $\cL$ is a collection of $q$ subsets of $\{1,\ldots,n\}$ and
each $m_i \in \{0, 1, \ldots, n\}$. In order for the triple to arise as the type of
an inequality it must satisfy a number of conditions given next:
\begin{itemize}
\item[(r.1)] $\cI \subseteq \{1, \ldots, q\}$. The sets $\{m_i\}\cup \cL_i$ ($i \in \cI$)are pairwise disjoint.  For any $i \in \cI$, if $|\cL_i| < q-1$ then $m_i \in \cL_i$, and if $m_i \notin \cL_i$ then $w_{m_i} < w_{m_k}$ for all $k > i$ with $k \in \cI$.
\item[ (r.2)] We require that $|\cL_i| < q^2$ for all $i \in \cI$.  In terms of separation from \textit{undominated, valid} inequalities requirement is valid in light of Theorem \ref{nodrag} and Lemma \ref{smallL}.
\end{itemize}

\noindent Stronger conditions can be imposed, however these assumptions suffice
to prove:
\begin{lemma} For given $q$, the set of pairs satisfying (r.1)-(r.2) includes all types arising from undominated, valid inequalities \eqref{genq}. The total number of tuples that satisfy (r.1)-(r.4) is at most $O(q^2 2^q n^{q^3})$.
\end{lemma}  
\noindent {\em Proof.} Follows from the above discussion and the fact that there
are at most $2^q$ choices for sets of indices $i$ with $m_i > 0$. \QED
\subsection{Separation using a given type}\label{sepgiventype}
Assume again a given $y \in [0,1]^n$. Consider a fixed triple $\tau = (\cI, \cL, m)$ that has been enumerated as indicated above, i.e, it satisfies (r.1)-(r.2). Here we will first describe an optimization problem whose solution either:
\begin{itemize}
\item [(a)] Proves that $y$ satisfies all inequalities \eqref{genq} of type $(\cI, \cL,m)$ that are valid for \eqref{minknap} (if any such inequalities exist), or
\item [(b)] Finds a valid inequality \eqref{genq}  for \eqref{minknap} that is violated by $y$.
\end{itemize}
To construct the formulation, we write, each $ i \in \cI$,
\begin{subequations}\label{condtau}\begin{align}
    M_i \ & \doteq \  \min \left\{ \, \min_{h \in \cL_i} w_h \ , \ \min_{k > i} w_{m_k} - 1 \, \right\}, \quad \cF_i \ \doteq \ \cup_{k \neq i} \{m_k\} \cup \cL_k\\
    \cR_i & \ \doteq \ \left\{\begin{tabular}{cl}
    $\left\{j \notin \cF_i \, : \, w_{m_i} \le w_j \le M_i \right\}$, & \mbox{if $|\cL_i| \ge q - 1$}\\
    & \\
    $\emptyset$, & \mbox{otherwise.} \end{tabular} \right.  \label{Ridef}\\
    \cV_i \ &  \doteq \  \cR_i \cup \cL_i.    \label{Videf}
\end{align}\end{subequations}

\noindent Now we describe our formulation. Let $\sigma(\tau)$ be the value of problem \eqref{sig} (computed exactly as in that formulation).
Then we solve the problem:
\begin{subequations} \label{pqtype} \begin{align}
    \Omega(y, \tau) \ &  \doteq  \ \min \quad \sum_{i \in \cI} i \left( \sum_{j \in \cV_{i}} y_j z_j \right) \\
    & \ \text{s.t.} \quad z_j = 1 \quad \forall \, j \in \cup_{i \in \cI} \{m_i\} \cup \cL_i \, \label{qfix1}\\
    & \ \qquad \ z_j = 0 \quad \forall \, j \notin \cup_i \cV_{i}, \label{qfix0}\\    
    &  \ \qquad  \sum_{j} w_j z_j \ \ge \ \sigma(\tau) \ + \ \sum_{j = 1}^n w_j - w_0 + 1, \label{thecovertype} \\    
    &  \ \qquad z \in \{0,1\}^n \nonumber
\end{align} \end{subequations}
\begin{lemma} \label{givenvalid} Suppose there is an inequality $\sum_{i = 1}^q i \, x(\cS_i) \ \ge \ q $
  of type $\tau$ valid for the knapsack \eqref{minknap}.  Then, setting $\hat z_j = 1$ iff $j \in \cup_i \cS_i$
  yields a feasible solution to \eqref{pqtype}.
\end{lemma}
\noindent {\em Proof.}  By definition of the type $\tau$ $\hat z$ satisfies \eqref{qfix1}.
  To show \eqref{qfix0} holds at $\hat z$, we will show that $\cS_i \subseteq \cV_i$ for all $i \in \cI$. Let $j \in \cS_i$;  if $j \in \{m_i\} \cup \cL_i$ then   by construction $j \in \cV_i$, so assume
$j \notin \{ m_i \} \cup \cL_i$. In this case we prove $j \in \cR_i$.  To do so, note that 
by definition
of $m_i$ in (t.2), $w_{m_i} \le w_j$. Moreover $j \in \cS_i \setminus \cL_i$ and hence (Remark \ref{booboo})
$w_j \le M_i$.  Finally, the $S_k$ are disjoint, and so $j \notin \cup_{k \neq i} \cS_k$ and hence $j \notin \cF_i$.  Thus indeed $j \in \cR_i$ as desired.
  To complete the proof we need to show that $\hat z$ satisfies \eqref{thecovertype}, but this follows
  from Lemma \ref{charqlemma}. \QED
  \begin{lemma} \label{validsogiven} Suppose 
    $\tilde z$ is a feasible solution for \eqref{pqtype}. 
    For $i \in \cI$ define $\cS_i \doteq \cV_i \cap \{j \, : \, \tilde z_j = 1\}$.  Then
\begin{align}     
  & \sum_{i \in \cI} i \, x(\cS_i) \ \ge \ q \label{proposed}
\end{align}  is valid for \eqref{minknap}.
  \end{lemma}
  \noindent {\em Proof.}  First note that for any $i \in \cI$, $\cS_i \supseteq \cL_i$ (by \eqref{qfix1}).  Moreover, for each $i \in \cI$,  either (a) $|\cL_i| < q - 1$ in which case $\cS_i = \cL_i$ or (b) $|\cL_i| \ge q - 1$ and $w_j \ge w_h$ for each $j \in \cL_i$ and $h \in \cS_i \setminus \cL_i$.  Thus it follows
  (Lemma \ref{restrictsig}) that if \eqref{proposed} has a certain type $\tau'$, then $\sigma(\tau') = \sigma(\tau)$.  As a result, constraint \eqref{thecovertype} and Lemma \ref{charqlemma} imply that \eqref{proposed} is valid for \eqref{minknap}. \QED

\noindent   We can now prove our key separation theorem.

  \begin{theorem} \label{separate}\hspace{.1in}
  The vector $y \in [0,1]^n$ violates an inequality of type $\tau$ valid for \eqref{minknap} iff  $\Omega(y, \tau) < q$.
\end{theorem}
\noindent {\em Proof.} Suppose first that $\sum_{i\in \cI}  i \, x(\cS_i)  \ \ge \ q $
is an inequality of type $\tau$, valid for \eqref{genq} and violated by $y$.  By Lemma \ref{givenvalid}, by  setting $\hat z_j = 1$ iff $j \in \cup_{i = 1}^q \cS_{e(i)}$ we obtain a feasible solution for problem \ref{pqtype}.  But since the objective value attained by $\hat z$ in this problem
equals
$\sum_{i \in \cI} i \, y(\cS_i) < q$ we conclude as desired. 

Now assume $\Omega(y, \tau) < q$. Let $\tilde z$ be an optimal solution for \eqref{pqtype}. By Lemma \ref{validsogiven} the inequality
\begin{align}     
  & \sum_{i \in \cI} i \, x(\cS_i) \ \ge \ q \label{proposed2}
\end{align}
is valid for \eqref{minknap} where for $i \in \cI$ we define $\cS_i \doteq \cV_i \cap \{j \, : \, \tilde z_j = 1\}$.
But since 
$$ q > \Omega(y, \tau) = \sum_{i \in \cI} i \, y(\cS_i)$$
we conclude $y$ violates \eqref{proposed2}. \QED

\subsubsection{Near separation in polynomial time}
In order to prove Theorem \ref{ourminknap} there remains the issue of the \textit{complexity} of solving problems of the form \eqref{pqtype}.  These are min-knapsack problems, for which an FPTAS
exists, based on that for the standard knapsack problem \cite{lawler}, \cite{ibarrakim}.  Relying
on such an FPTAS would yield a proof of Theorem \ref{ourminknap} (though,technically, the complexity would depend polynomially on $p/\epsilon$).  However this route would yield an algorithm that
relies on the traditional techniques: dynamic programming and coefficient scaling.

Here we indicate a simpler technique that applies in this case\footnote{We estimate that this is a folklore trick}. Consider, again, a given value of $q$ and an inquality of type $\tau$ as in the sections above. For $1 \le j \le n$ define $\hat y_j \doteq \frac{1}{q n^{2}} \lceil q n^{2} y_j \rceil$, i.e. the ``round-up'' of $y_j$ to the nearest integer multiple of
$\frac{1}{q n^{2}}$.  Then for any type $\tau$
$$ V(y, \tau) \le V(\hat y, \tau) \le V(y, \tau) + \frac{1}{n}$$
and so $V(\hat y , \tau) < q$ implies that $y$ violates an inequality of type
$\tau$ whereas if $V(\hat y , \tau) \ge q$ then $y$ satisfies every inequality
of type $\tau$ within additive error at most $1/n$, which is less than $\epsilon$ for $n$ large enough.

Moreover
$V(\hat y, \tau)$ can be computed in polynomial time, since it can be restated as a min-knapsack
problem with nonnegative, integral objective coefficients bounded above by $q n^{2}$.  Such
a min-knapsack problem can be solved using dynamic-programming (no need for coefficient scaling)\footnote{In fact even the dynamic-programming step can be eliminated \cite{bienben}.}. We have thus proved
Theorem \ref{ourminknap}.\\

\hspace{.1in}\\

\noindent {\bf Acknowledgement}. The work of the first author was partly funded by award ONR-GG012500.

\hspace {.1in}\\
\tiny Tue.Jun.19.151502.2018@blacknwhite

\normalsize
\bibliography{vb} 
\bibliographystyle{siam}

\end{document}